\documentclass[a4paper]{article}
\usepackage{makeidx}
\usepackage{amsfonts}
\usepackage{amsmath}
\usepackage{latexsym}
\usepackage{amssymb}
\usepackage{amsthm}
\usepackage{amscd}
\usepackage[dvips]{graphicx}
\usepackage[dvips]{color}

\title{Minor summation formula and a proof of Stanley's open problem}

\author{Masao Ishikawa\\
\small Faculty of Education, Tottori University\\[-0.8ex]
\small Koyama, Tottori, Japan\\[-0.8ex]
\small \texttt{ishikawa@fed.tottori-u.ac.jp}
}

\date{
}

\theoremstyle{definition}
\newtheorem{theorem}{Theorem}[section]
\newtheorem{prop}[theorem]{Proposition}
\newtheorem{lemma}[theorem]{Lemma}
\newtheorem{corollary}[theorem]{Corollary}

\newtheorem{remark}[theorem]{Remark}
\newtheorem{conjecture}[theorem]{Conjecture}
\newenvironment{demo}[1]{%
  \trivlist
  \item[\hskip\labelsep
        {\bf #1.}]
}{%
  \endtrivlist
}
\numberwithin{equation}{section} 

\begin{document}

\maketitle



                                
\newdimen\Squaresize \Squaresize=20pt
\newdimen\thickness \thickness=1pt         
                                                    
\def\Square#1{\hbox{\vrule width \thickness
   \vbox to \Squaresize{\hrule height \thickness\vss                            
      \hbox to \Squaresize{\hss#1\hss}
   \vss\hrule height\thickness} 
\unskip\vrule width \thickness} 
\kern-\thickness}                                                            
                               
\def\vsquare#1{\vbox{\Square{$#1$}}\kern-\thickness}
\def\blank{\omit\hskip\Squaresize}

\def\fibyoung#1{\let\\=\cr              
\vbox{\smallskip\offinterlineskip
\halign{&\vsquare{##}\cr #1}}\,}

%


\def\borderlessrect#1#2{\hbox{\hskip \thickness
   \vbox to \Squaresize{\vskip \thickness \vss
      \hbox to #2 {\hss #1\hss}
   \vss\vskip\thickness} 
\unskip\hskip \thickness} 
\kern-\thickness}                                                            
                               
\def\vborderlessrect#1#2{\vbox{\borderlessrect{$#1$}{#2}}\kern-\thickness}

\def\borderless#1{\omit\vborderlessrect{#1}{\Squaresize}}
\def\borderlessrc#1#2{\omit\vborderlessrect{#1}{#2}}

\def\msquare#1{\vbox{\hbox{\vrule width \thickness
   \vbox to \Squaresize{\hrule height \thickness
      \hbox to \Squaresize{\hfil{\sevenrm #1}}
   \vfil\hrule height\thickness}
\unskip\vrule width \thickness}
\kern-\thickness}\kern-\thickness}

\def\twosquare#1#2{\vbox{\hbox{\vrule width \thickness 
   \vbox to \Squaresize{\hrule height \thickness
      \hbox to \Squaresize{\hfil{\sevenrm #1}}\vss
      \hbox to \Squaresize{\hss{#2}\hss}
   \vfil\hrule height\thickness}
\unskip\vrule width \thickness}
\kern-\thickness}\kern-\thickness}

\def\twoblank#1#2{\vbox{\hbox{
   \vbox to \Squaresize{\vskip 2pt
      \hbox to \Squaresize{\hfil{\sevenrm #1}\ }\vss
      \hbox to \Squaresize{\hss{#2}\hss}
   \vfil}\unskip\kern-\thickness}
}\unskip\kern-\thickness}

\def\young#1{
\def\>{\blank}
\def\<{\borderless}
\def\*{\borderlessrc}
\def\p{\omit\msquare}
\def\t{\omit\twosquare}
\def\b{\omit\twoblank}
\let\\=\cr 
\vbox{\smallskip\offinterlineskip
\halign{&\vsquare{##}\cr #1}}}

\def\thinbox#1{\omit \thickness=.5pt \vsquare{#1}}
\def\tyoung{\def\t{\thinbox} \thickness=1pt \young}

\newdimen\smsquaresize \smsquaresize=12pt
\newdimen\smthickness \smthickness=.5pt
\font\smcellfont=cmss8 scaled \magstep0

\def\smsquare#1{\hbox{\vrule width \smthickness
   \unskip\vbox to \smsquaresize{\hrule height \smthickness\vss
      \hbox to \smsquaresize{\hss{\smcellfont #1}\hss}
   \vss\hrule height\smthickness} 
\unskip\vrule width \smthickness} 
\kern-\smthickness}

\def\smvsquare#1{\vbox{\smsquare{$#1$}}\kern-\smthickness}
\def\blank{\omit\hskip\smsquaresize}

\def\smborderlessrect#1#2{\hbox{\hskip \smthickness
   \vbox to \smsquaresize{\vskip \smthickness \vss
      \hbox to #2 {\hss #1\hss}
   \vss\vskip\smthickness} 
\unskip\hskip \smthickness} 
\kern-\smthickness}                                                            
                               
\def\smvborderlessrect#1#2{\vbox{\smborderlessrect{$#1$}{#2}}\kern-\smthickness}

\def\smborderless#1{\omit\smvborderlessrect{#1}{\smsquaresize}}
\def\smborderlessrc#1#2{\omit\smvborderlessrect{#1}{#2}}

\def\smyoung#1{
\def\<{\smborderless}
\def\*{\smborderlessrc}
\let\\=\cr 
\vbox{\smallskip\offinterlineskip
\halign{&\smvsquare{##}\cr #1}}}
\newdimen\vsmsquaresize \vsmsquaresize=10pt
\newdimen\vsmthickness \vsmthickness=.5pt
\font\vsmcellfont=cmsl8 scaled \magstep0
\font\vsmletterfont=cmr6 scaled \magstep0

\def\vsmsquare#1{\hbox{\vrule width \vsmthickness
   \unskip\vbox to \vsmsquaresize{\hrule height \vsmthickness\vss
      \hbox to \vsmsquaresize{\hss{\vsmcellfont #1}\hss}
   \vss\hrule height\vsmthickness} 
\unskip\vrule width \vsmthickness} 
\kern-\vsmthickness}
\def\vsmvsquare#1{\vbox{\vsmsquare{#1}}\kern-\vsmthickness}
\def\vsmblank{\omit\hskip\vsmsquaresize}
\def\vsmborderless#1{\hbox{\hskip \vsmthickness\unskip
   \vbox to \vsmsquaresize{\vss
      \hbox to \vsmsquaresize{\hss{\vsmletterfont #1}\hss}
   \vss} 
\unskip\hskip \vsmthickness} 
\kern-\vsmthickness}                                                            \def\vsmvborderless#1{\vbox{\vsmborderless{#1}}\kern-\vsmthickness}

\def\vsmyoung#1{
\def\>{\vsmblank}
\def\<{\omit\vsmvborderless}
\let\\=\cr 
\vbox{\smallskip\offinterlineskip
\halign{&\vsmvsquare{##}\cr #1}}}



\def\rdots{\mathinner{\mkern1mu\raise1pt\hbox{.}\mkern2mu\raise4pt\hbox{.}\mkern2mu\raise7pt\hbox{.}\mkern1mu}}

\def\defterm#1{{\sl #1}\/}

\def\module{\operatorname{mod}}

\def\ep{\varepsilon}
\def\lam{\lambda}

\def\Comp{\Bbb{C}}
\def\Nat{\Bbb{N}}
\def\Int{\Bbb{Z}}
\def\Pos{\Bbb{P}}

\def\Pf{\operatorname{Pf}}
\def\pf{\operatorname{pf}}
\def\odd{\operatorname{odd}}
\def\even{\operatorname{even}}
\def\sgn{\operatorname{sgn}}
\def\MOD{{\,\operatorname{mod}\,}}
\def\Res{\operatorname{Res}}

\def\half{\frac12}
\def\trans{{}^t\!}
\def\cc#1{{\ooalign{\hfil\raise-.02ex\hbox{#1}\hfil\crcr\mathhexbox20D}}}
\def\cl{\mathhexbox20D}
\def\cb{\bullet}

\def\diag{\operatorname{diag}}

\def\Ind#1{\Cal{I}_{#1}}





\def\id{\operatorname{id}}
\def\wt{\operatorname{wt}}

          





\abstract{
In the open problem session of the FPSAC'03,
R.P.~Stanley gave an open problem about a certain sum of the Schur functions
(See \cite{Sta3}).
The purpose of this paper is to give a proof of this open problem.
The proof consists of three steps.
At the first step
we express the sum by a Pfaffian as an application of our minor summation formula
(\cite{IW1}).
In the second step we prove a Pfaffian analogue of Cauchy type identity
which generalize \cite{Su2}.
Then we give a proof of Stanley's open problem in Section~4.
At the end of this paper we present certain corollaries obtained from this identity
involving the Big Schur functions and some polynomials arising from the Macdonald polynomials,
which generalize Stanley's open problem.

}

\bigbreak

%
%
\section{Introduction}\label{sec:intro}

In the open problem session of the 15th Anniversary International Conference on Formal Power Series and Algebraic Combinatorics
(Vadstena, Sweden, 25 June 2003),
R.P.~Stanley gave an open problem on a sum of Schur functions 
with a weight including four parameters,
i.e. Theorem~\ref{conj:general}
(See \cite{Sta3}).
The purpose of this paper is to give a proof of this open problem.
In the process of our proof, we obtain a Pfaffian identity,
i.e. Theorem~\ref{fundamental},
which generalize the Pfaffian identities in \cite{Su2}.
Note that certain determinant and Pfaffian identities of this type first appeared in \cite{O2},
and applied to solve some alternating sign matrices enumerations under certain symmetries stated in \cite{Ku1}.
Certain conjectures which intensively generalize the determinant and Pfaffian identities of this type were stated in \cite{O4},
and a proof of the conjectured determinant and Pfaffian identities was given in \cite{IOTZ}.
Now we know that various methods may be adopted to prove this type of identity.
We can prove it algebraically using Dodgson's formula or the usual expansion formula of Pfaffians.
Here we state an analytic proof since this proof is due to the author and is not stated in other places.
Our proof proceeds by three steps.
In the first step we utilize the minor summation formula (\cite{IW1}) to express the sum of Schur functions as a Pfaffian.
In the second step we express the Pfaffian by a determinant using a Cauchy type Pfaffian formula (also see \cite{O3}, \cite{O4} and \cite{IOTZ}),
and try to simplify it as much as possible.
In the process of this step,
it is conceivable that the determinants we treat may be closely related to characters of representations of $\operatorname{SP}_{2n}$ and $\operatorname{SO}_{m}$
(See \cite{FH}, \cite{IOW} and \cite{IW5}).
In the final step we complete our proof using a key proposition,
i.e. Proposition~\ref{key_lemma}
(See \cite{Sta1} and \cite{Ste2}).
At the end of this paper
we state some corollaries which generalize Stanley's open problem to the big Schur functions,
and to certain polynomials arising from the Macdonald polynomials.
Furthermore,
in the forthcoming paper \cite{IZ},
we study a finite version of Boulet's theorem and present certain relations with
orthogonal polynomials and the basic hypergeometric series.
In the paper we find more applications of the Pfaffain expression of Stanley's weight $\omega(\lambda)$
obtained in this paper,
and also study a certain summation of Schur's $Q$-functions weighted by $\omega(\lambda)$.

We follow the notation in \cite{Ma} concerning symmetric functions.
In this paper
we use a symmetric function $f$ in $n$ variables $(x_1,\dots,x_n)$,
which is usually written as $f(x_1,\dots,x_n)$,
and also a symmetric function $f$ in countably many variables $x=(x_1,x_2,\dots)$,
which is written as $f(x)$
(for detailed description of the ring of symmetric functions in countably many variables,
see \cite{Ma}, I, sec.2).
To simplify this notation
we express the $n$-tuple $(x_1,\dots,x_n)$ by $X_n$,
and sometimes simply write $f(X_n)$ for $f(x_1,\dots,x_n)$.
When the number of variables is finite 
and there is no fear of confusion what this number is,
we simply write $X$ for $X_n$ in abbreviation.
Thus $f(x)$ is in countably many variables,
but $f(X)$ is in finitely many variables and the number of variables is clear from the assumption.

Given a partition $\lambda$,
define $\omega(\lambda)$ by
\begin{equation*}
\omega(\lambda)=a^{\sum_{i\geq1}\lceil\lambda_{2i-1}/2\rceil}
b^{\sum_{i\geq1}\lfloor\lambda_{2i-1}/2\rfloor}
c^{\sum_{i\geq1}\lceil\lambda_{2i}/2\rceil}
d^{\sum_{i\geq1}\lfloor\lambda_{2i}/2\rfloor},
\end{equation*}
where $a$, $b$, $c$ and $d$ are indeterminates,
 and $\lceil x\rceil$ (resp. $\lfloor x\rfloor$) stands for the smallest (resp. largest) integer greater (resp. less) than or equal to $x$ 
for a given real number $x$.
For example,
if $\lambda=(5,4,4,1)$
then
$\omega(\lambda)$ is the product of the entries in the following diagram for $\lambda$.
\[
\young{
a&b&a&b&a\\
c&d&c&d\\
a&b&a&b\\
c\\
}
\]
Let $s_{\lambda}(x)$ denote the Schur function corresponding to a partition $\lambda$.
R.~P.~Stanley gave the following conjecture in the open problem session of FPSAC'03.
\begin{theorem}
Let
\begin{equation*}
z=\sum_{\lambda}\omega(\lambda)s_{\lambda}.
\end{equation*}
Here the sum runs over all partitions $\lambda$.
Then we have
\begin{align}
\log z-\sum_{n\geq1}\frac1{2n}a^n(b^n-c^n)p_{2n}
-\sum_{n\geq1}\frac1{4n}a^nb^nc^nd^np_{2n}^2\nonumber\\
\in\Bbb{Q}[[p_1,p_3,p_5,\dots]].
\label{conj:general}
\end{align}
Here $p_r=\sum_{i\geq1} x_i^r$ denote the $r$th power sum symmetric function.
\end{theorem}
As a special case of this open problem,
if we put $b=c=a^{-1}$ and $d=a$,
and check the constant term of the both sides,
then we obtain the following simple case:
\begin{corollary}
\label{conj:simple}
Let
\begin{equation*}
y=\sum_{{\lambda}\atop{\lambda,\lambda'\text{ even}}}s_{\lambda}(x).
\end{equation*}
Here the sum runs over all partitions $\lambda$ such that $\lambda$ and $\lambda'$ are even partitions (i.e. with all parts even).
Then we have
\begin{equation}
\log y-\sum_{n\geq1}\frac1{4n}p_{2n}^2\in \Bbb{Q}[[p_1,p_3,p_5,\dots]].
\end{equation}
\end{corollary}

In the rest of this section we briefly recall the definition of Pfaffians.
For a detailed explanation of Pfaffians,
the reader can consult \cite{Kn} and \cite{Ste1}.
Let $n$ be a non-negative integer
and assume we are given a $2n$ by $2n$ skew-symmetric matrix $A=(a_{ij})_{1\le i,j\le 2n}$, (i.e. $a_{ji}=-a_{ij}$),
whose entries $a_{ij}$ are in a commutative ring.
The \defterm{Pfaffian} of $A$ is,
by definition,
\begin{equation*}
\label{def_pfaffian}
\Pf(A)=\sum \epsilon(\sigma_{1},\sigma_{2},\hdots,\sigma_{2n-1},\sigma_{2n})\,
a_{\sigma_{1}\sigma_{2}} \dots a_{\sigma_{2n-1}\sigma_{2n}}.
\end{equation*}
where the summation is over all partitions $\{\{\sigma_{1},\sigma_{2}\}_{<},\hdots,\{\sigma_{2n-1},\sigma_{2n}\}_{<}\}$
of $[2n]$ into $2$-elements blocks,
and where $\epsilon(\sigma_{1},\sigma_{2},\hdots,\sigma_{2n-1},\sigma_{2n})$ denotes the sign of the permutation
\begin{equation*}
\begin{pmatrix}
1&2&\cdots&2n\\
\sigma_{1}&\sigma_{2}&\cdots&\sigma_{2n}
\end{pmatrix}.
\end{equation*}
We call a partition $\{\{\sigma_{1},\sigma_{2}\}_{<},\hdots,\{\sigma_{2n-1},\sigma_{2n}\}_{<}\}$
of $[2n]$ into $2$-elements blocks \defterm{a matching} or \defterm{$1$-factor} of $[2n]$.

%
%
\section{Minor Summation Formula}\label{sec:msf}

First we restrict our attention to the finite variables case.
Let $n$ be a non-negative integer.
We put
\begin{equation}
y_n=y_n(X_{2n})
=\sum_{{\lambda}\atop{\lambda,\lambda'\text{ even}}}s_{\lambda}(X_{2n})
=\sum_{{\lambda}\atop{\lambda,\lambda'\text{ even}}}s_{\lambda}(x_1,\dots,x_{2n}).
\end{equation}
where $s_{\lambda}(X_{2n})$ is the Schur function corresponding to a partition $\lambda$ in $2n$ variables $x_1,\dots,x_{2n}$.
Then there is a known formula which is originally due to Littlewood as follows.
(See \cite{Su2}).
\begin{align}
y_n(X_{2n})
=\frac1{\prod_{1\leq i<j\leq 2n}(x_i-x_j)}\Pf\left[\frac{x_i-x_j}{1-x_i^2x_j^2}\right]_{1\leq i<j\leq 2n}.
\end{align}
The aim of this section is to prove the following theorem which generalize this identity.

\begin{theorem}
\label{theorem:pfaffian}
Let $n$ be a positive integer and let $\omega(\lambda)$ be as defined in Section~\ref{sec:intro}.
Let
\begin{equation}
z_n
=z_n(X_{2n})
=\sum_{\ell(\lambda)\leq 2n}\omega(\lambda)s_{\lambda}(X_{2n})
=\sum_{\ell(\lambda)\leq 2n}\omega(\lambda)s_{\lambda}(x_{1},\dots,x_{2n})
\end{equation}
be the sum restricted to $2n$ variables.
Then we have
\begin{align}
\label{eq:general}
z_{n}(X_{2n})=\frac{(abcd)^{-\binom{n}2}}{\prod_{1\leq i<j\leq2n}(x_i-x_j)}
\Pf\left(p_{ij}\right)_{1\leq i<j\leq 2n},
\end{align}
where $p_{ij}$ is defined by
\begin{align}
\label{eq:entries}
p_{ij}
=\frac{\begin{vmatrix}
x_{i}+a x_{i}^2&1-a(b+c)x_{i}-abc x_{i}^3\\
x_{j}+a x_{j}^2&1-a(b+c)x_{j}-abc x_{j}^3
\end{vmatrix}}
{(1-abx_i^2)(1-abx_j^2)(1-abcdx_i^2x_j^2)}.
\end{align}
\end{theorem}

Let $m$, $n$ and $r$ be integers such that $r\leq m,n$ and let $T$ be an $m$ by $n$ matrix.
For any index sets $I=\{i_1,\dots,i_r \}_{<}\subseteq[m]$ and $J=\{j_1,\dots,j_r\}_{<}\subseteq[n]$,
let $\Delta^{I}_{J}(A)$ denote the submatrix obtained by selecting the rows indexed by $I$ and the columns indexed by $J$.
If $r=m$ and $I=[m]$,
we simply write $\Delta_{J}(A)$ for $\Delta^{[m]}_{J}(A)$.
Similarly,
if $r=n$ and $J=[n]$,
we write $\Delta^{I}(A)$ for $\Delta_{[n]}^{I}(A)$.
For any finite set $S$ and a non-negative integer $r$,
let $\binom{S}{r}$ denote the set of all $r$-element subsets of $S$.
We cite a theorem from \cite{IW1} which we call a minor summation formula:
\begin{theorem}
\label{msf}
Let $n$ and $N$ be non-negative integers such that $2n\le N$.
Let $T=(t_{ij})_{1\leq i\leq2n, 1\leq j\leq N}$ be a $2n$ by $N$ rectangular matrix,
and let $A=(a_{ij})_{1\le i,j\le N}$ be a skew-symmetric matrix of size $N$.
Then
\begin{align*}
\label{eq_msf}
\sum_{I\in\binom{[N]}{2n}}
\Pf\left(\Delta^{I}_{I}(A)\right) \det\left(\Delta_I(T)\right)&
=\Pf\left(TA\,{}^t\kern-1pt T\right).
\end{align*}
If we put $Q=\left(Q_{ij}\right)_{1\leq i,j\leq 2n}=TA\,{}^t\kern-1pt T$, then its entries are given by
\begin{equation*}
Q_{ij}=\sum_{1\le k<l\le N} a_{kl} \det\left(\Delta^{ij}_{kl}(T)\right),
\qquad(1\le i,j\le 2n).
\end{equation*}
Here we write $\Delta^{ij}_{kl}(T)$ for $\Delta^{\{ij\}}_{\{kl\}}(T)=\begin{vmatrix}t_{ik}&t_{il}\\t_{jk}&t_{jl}\end{vmatrix}$.
$\Box$
\end{theorem}

Before we proceed to the proof of Theorem~\ref{theorem:pfaffian},
we cite a lemma from \cite{IW1}.
The proof is not difficult,
but we omit the proof and the reader can consult \cite{IW1}, Section~4, Lemma~7.
\begin{lemma}
\label{lemma:product}
Let $x_i$ and $y_j$ be indeterminates,
and let $n$ be a non-negative integer.
Then
\begin{equation}
\label{prod_pf}
\Pf[x_iy_j]_{1\leq i<j\leq 2n}=\prod_{i=1}^{n}x_{2i-1}\prod_{i=1}^{n}y_{2i}.
\ \Box
\end{equation}
\end{lemma}

Given a partition $\lambda=(\lambda_1,\dots,\lambda_{m})$ satisfying $\ell(\lambda)\leq m$,
we associate a decreasing sequence $\lambda+\delta_m$ which is usually denoted by $l=(l_1,\dots,l_m)$,
where $\delta_m=(m-1,m-2,\dots,0)$.

\begin{lemma}
\label{coefficients}
Let $n$ be a non-negative integer.
Let $\lambda=(\lambda_1,\dots,\lambda_{2n})$ be a partition such that $\ell(\lambda)\leq2n$,
and put $l=(l_1,\dots,l_{2n})=\lambda+\delta_{2n}$.
Define a $2n$ by $2n$ skew-symmetric matrix $A=(\alpha_{ij})_{1\leq i,j\leq2n}$
by
\begin{equation*}
\alpha_{ij}=a^{\lceil(l_i-1)/2\rceil}b^{\lfloor(l_i-1)/2\rfloor}c^{\lceil l_j/2\rceil}d^{\lfloor l_j/2\rfloor}
\end{equation*}
for $i<j$,
and as $\alpha_{ji}=-\alpha_{ij}$ holds for any $1\leq i,j\leq2n$.
Then we have
\begin{equation*}
\Pf\left[A\right]_{1\leq i,j\leq 2n}
=(abcd)^{\binom{n}2}\omega(\lambda).
\end{equation*}
\end{lemma}
\begin{demo}{Proof}
By Lemma~\ref{lemma:product},
we have
\begin{equation*}
\Pf[A]
=\prod_{i=1}^{n}a^{\lceil(l_{2i-1}-1)/2\rceil}b^{\lfloor(l_{2i-1}-1)/2\rfloor}
\prod_{j=1}^{n}c^{\lceil l_{2j}/2\rceil}d^{\lfloor l_{2j}/2\rfloor}.
\end{equation*}
Since $l_{2i-1}-1=\lambda_{2i-1}+2(n-i)$ and $l_{2j}=\lambda_{2j}+2(n-j)$,
this Pfaffian becomes
\begin{equation*}
\prod_{i=1}^{n}a^{\lceil\lambda_{2i-1}/2\rceil+n-i}b^{\lfloor\lambda_{2i-1}/2\rfloor+n-i}
\prod_{j=1}^{n}c^{\lceil\lambda_{2j}/2\rceil+n-j}d^{\lfloor\lambda_{2j}/2\rfloor+n-j},
\end{equation*}
which is easily seen to be $(abcd)^{\binom{n}2}\omega(\lambda)$.
$\Box$
\end{demo}
Now we are in the position to give a proof of Theorem~\ref{theorem:pfaffian}.
\begin{demo}{Proof of Theorem~\ref{theorem:pfaffian}}
By Theorem~\ref{msf}
it is enough to compute
\[
\beta_{ij}=\sum_{k\geq l\geq0}
a^{\lceil(k-1)/2\rceil}b^{\lfloor(k-1)/2\rfloor}
c^{\lceil l/2\rceil}d^{\lfloor l/2\rfloor}
\begin{vmatrix}
x_{i}^{k}&x_{i}^{l}\\
x_{j}^{k}&x_{j}^{l}
\end{vmatrix}.
\]
Let
$
f^{ij}_{kl}=
a^{\lceil(k-1)/2\rceil}b^{\lfloor(k-1)/2\rfloor}
c^{\lceil l/2\rceil}d^{\lfloor l/2\rfloor}
\begin{vmatrix}
x_{i}^{k}&x_{i}^{l}\\
x_{j}^{k}&x_{j}^{l}
\end{vmatrix},
$
then, this sum can be divided into four cases, i.e.
\[
\beta_{ij}=
\sum_{{k=2r+1,\,l=2s}\atop{r\geq s\geq0}}f^{ij}_{kl}
+\sum_{{k=2r,\,l=2s}\atop{r\geq s\geq0}}f^{ij}_{kl}
+\sum_{{k=2r+1,\,l=2s+1}\atop{r\geq s\geq0}}f^{ij}_{kl}
+\sum_{{k=2r+2,\,l=2s+1}\atop{r\geq s\geq0}}f^{ij}_{kl}.
\]
We compute each case:
\begin{enumerate}
\item[(i)]
If $k=2r+1$ and $l=2s$ for $r\geq s\geq0$,
then
\begin{eqnarray*}
\sum_{{k=2r+1,\,l=2s}\atop{r\geq s\geq0}}f^{ij}_{kl}
&&=\sum_{r\geq s\geq0}
a^rb^rc^sd^s
\begin{vmatrix}
x_{i}^{2r+1}&x_{i}^{2s}\\
x_{j}^{2r+1}&x_{j}^{2s}
\end{vmatrix}\\
&&=\sum_{r\geq s\geq0}
c^sd^s
\begin{vmatrix}
\frac{a^sb^sx_{i}^{2s+1}}{1-abx_{i}^2}&x_{i}^{2s}\\
\frac{a^sb^sx_{j}^{2s+1}}{1-abx_{j}^2}&x_{j}^{2s}
\end{vmatrix}\\
&&=\frac{(x_{i}-x_{j})(1+abx_{i}x_{j})}{(1-abx_{i}^2)(1-abx_{j}^2)(1-abcdx_{i}^2x_{j}^2)}.
\end{eqnarray*}
In the same way we obtain the followings by straight forward computations.
\item[(ii)]
If $k=2r$ and $l=2s$ for $r\geq s\geq0$,
then
\begin{eqnarray*}
\sum_{{k=2r,\,l=2s}\atop{r\geq s\geq0}}f^{ij}_{kl}
&&=\frac{a(x_{i}^2-x_{j}^2)}{(1-abx_{i}^2)(1-abx_{j}^2)(1-abcdx_{i}^2x_{j}^2)}.
\end{eqnarray*}
\item[(iii)]
If $k=2r+1$ and $l=2s+1$ for $r\geq s\geq0$,
then
\begin{eqnarray*}
\sum_{{k=2r+1,\,l=2s+1}\atop{r\geq s\geq0}}f^{ij}_{kl}
&&=\frac{abc x_{i}x_{j}(x_{i}^2-x_{j}^2)}{(1-abx_{i}^2)(1-abx_{j}^2)(1-abcdx_{i}^2x_{j}^2)}.
\end{eqnarray*}
\item[(iv)]
If $k=2r+2$ and $l=2s+1$ for $r\geq s\geq0$,
then
\begin{eqnarray*}
\sum_{{k=2r+2,\,l=2s+1}\atop{r\geq s\geq0}}f^{ij}_{kl}
&&=\frac{ac x_{i}x_{j}(x_{i}-x_{j})(1+abx_{i}x_{j})}{(1-abx_{i}^2)(1-abx_{j}^2)(1-abcdx_{i}^2x_{j}^2)}.
\end{eqnarray*}
\end{enumerate}
Summing up these four identities,
we obtain
\begin{align*}
\beta_{ij}
=\frac{(x_i-x_j)\{1+abx_ix_j+a(x_i+x_j)+abcx_ix_j(x_i+x_j)+acx_ix_j(1+abx_ix_j)\}}{(1-abx_i^2)(1-abx_j^2)(1-abcdx_i^2x_j^2)}.
\end{align*}
It is easy to see the numerator is written by the determinant,
and this completes the proof.
$\Box$
\end{demo}

%
%
\section{Cauchy Type Pfaffian Formulas}
The aim of this section is to derive \thetag{\ref{eq:cauchy}}.
In the next section we will use this identity to prove Stanley's open problem.
First we prove a fundamental Pfaffian identity,
i.e.
Theorem~\ref{fundamental},
and deduce all the identities in this section from this theorem.
In the latter half of this section
we also show that we can derive Sundquist's Pfaffian identities obtained in \cite{Su2}
from our theorem
although these identities have no direct relation to Stanley's open problem.
In that sense our theorem can be regarded as a generalization of Sundquist's Pfaffian identities (See \cite{Su2}).
An intensive generalization was conjectured in \cite{O4} and proved in \cite{IOTZ}.
There must be several ways to prove this type of identity.
In \cite{Su2} Sundquist gave a combinatorial proof of his Pfaffian identities.
In \cite{IOTZ} the authors adopted an algebraic method to prove identities of this type.
Here we give an analytic proof of our theorem,
in which we regard both sides of this identity as meromorphic functions
and check the Laurent series expansion at each isolated pole in the Riemann sphere.
The idea to use complex analysis to prove various determinant and Pfaffian identities is first hinted by Prof. H.~Kawamuko
to the author,
and the author recognized this can be a powerful tool to prove various identities including determinants and Pfaffians.
This idea was also used to prove a Pfaffian-Hafnian analogue of Borchardt's identity in \cite{IKO}.
In this section,
we first state our theorems and later give proofs of them.

First we fix notation.
Let $n$ be an non-negative integer.
Let $X=(x_{1},\dots,x_{2n})$, $Y=(y_{1},\dots,y_{2n})$, $A=(a_{1},\dots,a_{2n})$ and $B=(b_{1},\dots,b_{2n})$ be $2n$-tuples of variables.
Set $V^{n}_{ij}(X,Y;A,B)$ to be
\begin{equation*}
\begin{cases}
a_{i}x_{i}^{n-j}y_{i}^{j-1}
&\text{ if $1\leq j\leq n$,}\\
b_{i}x_{i}^{2n-j}y_{i}^{j-n-1}
&\text{ if $n+1\leq j\leq 2n$,}
\end{cases}
\end{equation*}
for $1\leq i\leq2n$,
and define $V^{n}(X,Y;A,B)$ by
\begin{equation*}
V^{n}(X,Y;A,B)=
\det\left(V^{n}_{ij}(X,Y;A,B)\right)_{1\leq i,j\leq2n}.
\end{equation*}
For example,
if $n=1$, then we have
$V^{1}(X,Y;A,B)=\begin{vmatrix}a_1&b_1\\a_2&b_2\end{vmatrix}$,
and
if $n=2$, then $V^{2}(X,Y;A,B)$ looks as follows:
\begin{eqnarray*}
V^{2}(X,Y;A,B)=\begin{vmatrix}
a_{1}x_{1}&a_{1}y_{1}&
b_{1}x_{1}&b_{1}y_{1}\\
a_{2}x_{2}&a_{2}y_{2}&
b_{2}x_{2}&b_{2}y_{2}\\
a_{3}x_{3}&a_{3}y_{3}&
b_{3}x_{3}&b_{3}y_{3}\\
a_{4}x_{4}&a_{4}y_{4}&
b_{4}x_{4}&b_{4}y_{4}\\
\end{vmatrix}.
\end{eqnarray*}

The main result of this section is the following theorem.
\begin{theorem}
\label{fundamental}
Let $n$ be a positive integer.
Let $X=(x_{1},\dots,x_{2n})$, $Y=(y_{1},\dots,y_{2n})$, $A=(a_{1},\dots,a_{2n})$, $B=(b_{1},\dots,b_{2n})$, $C=(c_{1},\dots,c_{2n})$ and $D=(d_{1},\dots,d_{2n})$
be $2n$-tuples of variables.
Then
\begin{align}
\label{eq:fundamental}
&\Pf\left[\frac{\begin{vmatrix}a_{i}&b_{i}\\a_{j}&b_{j}\end{vmatrix}
\cdot\begin{vmatrix}c_{i}&d_{i}\\c_{j}&d_{j}\end{vmatrix}}
{\begin{vmatrix}x_{i}&y_{i}\\x_{j}&y_{j}\end{vmatrix}}\right]_{1\leq i<j\leq 2n}
=\frac{V^{n}(X,Y;A,B)V^{n}(X,Y;C,D)}
{\displaystyle\prod_{1\leq i<j\leq 2n}\begin{vmatrix}x_{i}&y_{i}\\x_{j}&y_{j}\end{vmatrix}}.
\end{align}
\end{theorem}
The following proposition is obtained easily by elementary transformations of the matrices
and we will prove it later.
\begin{prop}
\label{prop:subs}
Let $n$ be a positive integer.
Let $X=(x_{1},\dots,x_{2n})$ be a $2n$-tuple of variables and let $t$ be an indeterminate.
Then
\begin{equation}
\label{eq:subs}
V^{n}(X,\pmb{1}+tX^2;X,\pmb{1})
=(-1)^{\binom{n}2}t^{\binom{n}2}\prod_{1\leq i<j\leq2n}(x_{i}-x_{j}),
\end{equation}
where $\pmb1$ denotes the $2n$-tuple $(1,\dots,1)$,
and $\pmb{1}+tX^2$ denotes the $2n$-tuple $(1+tx_{1}^2,\dots,1+tx_{2n}^2)$.
\end{prop}

Let $t$ be an arbitrary indeterminate.
If we set $y_{i}=1+tx_{i}^2$ in \thetag{\ref{eq:fundamental}}, then
\[
\begin{vmatrix}x_{i}&1+tx_{i}^2\\x_{j}&1+tx_{j}^2\end{vmatrix}
=(x_{i}-x_{j})(1-tx_{i}x_{j})
\]
and \thetag{\ref{eq:subs}} immediately implies the following corollary.
\begin{corollary}
\label{cor:pfaff-det}
Let $n$ be a non-negative integer.
Let $X=(x_1,\dots,x_{2n})$, $A=(a_1,\dots,a_{2n})$, $B=(b_1,\dots,b_{2n})$, $C=(c_1,\dots,c_{2n})$ and $D=(d_1,\dots,d_{2n})$
be $2n$-tuples of variables.
Then
\begin{align}
&\Pf\left[\frac{(a_ib_j-a_jb_i)(c_id_j-c_jd_i)}
{(x_i-x_j)(1-tx_ix_j)}\right]_{1\leq i<j\leq 2n}
\nonumber\\
&\qquad\qquad=\frac{V^{n}(X,\pmb{1}+tX^2;A,B)V^{n}(X,\pmb{1}+tX^2;C,D)}
{\prod_{1\leq i<j\leq 2n}(x_i-x_j)(1-tx_ix_j)}.
\end{align}
In particular, we have
\begin{align}
\label{eq:cauchy}
&\Pf\left[\frac{a_ib_j-a_jb_i}
{1-tx_ix_j}\right]_{1\leq i<j\leq 2n}
=(-1)^{\binom{n}2}t^{\binom{n}2}\frac{V^{n}(X,\pmb{1}+tX^2;A,B)}{\prod_{1\leq i<j\leq 2n}(1-tx_ix_j)}.
\ \Box
\end{align}
\end{corollary}

In the latter half of this section we show that we can derive Sundquist's Pfaffian identities from ours.
Let $X=(x_{1},\dots,x_{2n})$ and $A=(a_{1},\dots,a_{2n})$ be $2n$-tuples of variables
and let $S_{2n}$ act on each by permuting indices.
For compositions $\alpha$ and $\beta$ of length $n$,
we write
\begin{equation*}
a_{\alpha,\beta}(X;A)
=\sum_{\sigma\in S_{2n}}\epsilon(\sigma)\sigma\left(
a_{1} x_1^{\alpha_{1}}\cdots a_{n} x_{n}^{\alpha_{n}}
x_{n+1}^{\beta_{1}}\cdots x_{2n}^{\beta_{n}}
\right).
\end{equation*}
This is to say
\begin{equation*}
a_{\alpha,\beta}(X;A)
=\begin{vmatrix}
a_{1} x_{1}^{\alpha_1}& \hdots & a_{1} x_{1}^{\alpha_n} &
x_{1}^{\beta_1} & \hdots & x_{1}^{\beta_n}\\
\vdots & \ddots & \vdots & \vdots & \ddots & \vdots\\
a_{2n} x_{2n}^{\alpha_1}& \hdots & a_{2n} x_{2n}^{\alpha_n} &
x_{2n}^{\beta_1} & \hdots &  x_{2n}^{\beta_n}\\
\end{vmatrix}.
\end{equation*}
Let $n$ be a non-negative integer.
Let $\mathcal{P}_n$ denote the set of all partitions of the form
$\lambda=(\alpha_1,\dots,\alpha_r|\alpha_1+1,\dots,\alpha_r+1)$ in Frobenius notation with $\alpha_r\leq n-1$.
For example,
\begin{equation*}
\mathcal{P}_4=\{
\emptyset,1^2,21^2,2^3,31^3,32^21,3^22^2,3^4
\}.
\end{equation*}
We put
\begin{equation*}
U^{n}(X;A)=\sum_{\lambda\in\mathcal{P}_{n}}\sum_{\mu\in\mathcal{P}_{n}}
a_{\lambda,\mu}(X;A).
\end{equation*}
\begin{theorem}
\label{det_exp}
Let $n$ be a non-negative integer.
Let $X=(x_1,\dots,x_{2n})$ and $A=(a_1,\dots,a_{2n})$ be $2n$-tuples of variables.
Then
\begin{eqnarray}
U^{n}(X;A)=V^{n}(X,\pmb{1}-X^2;A,\pmb{1}),
\end{eqnarray}
where $\pmb1-X^2=(1-x_{1}^2,\dots,1-x_{2n}^2)$ and $\pmb1=(1,\dots,1)$.
\end{theorem}
From Theorem~\ref{fundamental}, Corollary~\ref{cor:pfaff-det} and Theorem~\ref{det_exp}
we obtain the following Pfaffian identities which are obtained in \cite{Su2}
(See \cite{Su2}, Theorem~2.1).
\begin{corollary}
(Sundquist)
\label{cor:Sunduist}
\begin{align}
\label{Sundquist1}
&\Pf\left[\frac{a_i-a_j}{x_i+x_j}\right]_{1\leq i<j\leq 2n}
=(-1)^{\binom{n}2}
\frac{a_{2\delta_n,2\delta_n}(X;A)}
{\prod_{1\leq i<j\leq 2n}(x_i+x_j)}
\\
\label{Sundquist2}
&\Pf\left[\frac{a_i-a_j}{1+x_ix_j}\right]_{1\leq i<j\leq 2n}
=\frac{U^{n}(X;A)}
{\prod_{1\leq i<j\leq 2n}(1+x_ix_j)}.
\end{align}
\end{corollary}

Now we state the proofs of our theorems.
Before we prove Theorem~\ref{fundamental},
we need two lemmas.
Let $n$ and $r$ be integers such that $2n\geq r\geq0$.
Let $X=(x_{1},\dots,x_{2n})$ be a $2n$-tuple of variables
and let $1\leq k_1<\dots<k_r\leq2n$ be a sequence of integers.
Let $X^{(k_1,\dots,k_r)}$ denote the $(2n-r)$-tuple of variables
obtained by removing the variables $x_{k_1}$, $\dots$, $x_{k_r}$ from $X_{2n}$.
\begin{lemma}
\label{substitution}
Let $n$ be a positive integer.
Let $X=(x_{1},\dots,x_{2n})$ and $A=(a_1,\dots,a_{2n})$
be $2n$-tuples of variables.
Let $k$, $l$ be any integers such that $1\leq k<l\leq2n$.
Then
\begin{align*}
V^{n}(X,\pmb{1};A,\pmb{1})\Big\rvert_{x_{l}=x_{k}}
&=(-1)^{k+l+n}(a_{k}-a_{l})\\
&\times
\prod_{\scriptstyle i=1\atop\scriptstyle i\neq k,l}^{2n}
(x_{i}-x_{k})
\cdot V^{n-1}(X^{(k,l)},\pmb{1}^{(k,l)};A^{(k,l)},\pmb{1}^{(k,l)}),
\end{align*}
where $\pmb{1}$ denotes the $2n$-tuple $(1,\dots,1)$.
\end{lemma}
\begin{demo}{Proof}
Without loss of generality we may assume that  $k=2n-1$ and $l=2n$.
From the definition,
$V^{n}(X,\pmb{1};A,\pmb{1})$ is in the form of
\begin{eqnarray*}
&&\det\left(
\begin{cases}
a_{i}x_{i}^{n-j}
&\text{ if $1\leq j\leq n$,}\\
x_{i}^{2n-j}
&\text{ if $n+1\leq j\leq 2n$.}\\
\end{cases}
\right)_{1\leq i,j\leq2n}.
\end{eqnarray*}
For example,
when $n=3$,
if we substitute $x_{6}=x_{5}$ into this determinant,
we obtain
\begin{eqnarray*}
V^{n}(X,\pmb{1};A,\pmb{1})\Big\rvert_{x_{6}=x_{5}}=
\begin{vmatrix}
a_{1}x_{1}^{2}&a_{1}x_{1}&a_{1}&
x_{1}^{2}&x_{1}&1\\
a_{2}x_{2}^{2}&a_{2}x_{2}&a_{2}&
x_{2}^{2}&x_{2}&1\\
a_{3}x_{3}^{2}&a_{3}x_{3}&a_{3}&
x_{3}^{2}&x_{3}&1\\
a_{4}x_{4}^{2}&a_{4}x_{4}&a_{4}&
x_{4}^{2}&x_{4}&1\\
a_{5}x_{5}^{2}&a_{5}x_{5}&a_{5}&
x_{5}^{2}&x_{5}&1\\
a_{6}x_{5}^{2}&a_{6}x_{5}&a_{6}&
x_{5}^{2}&x_{5}&1\\
\end{vmatrix}.
\end{eqnarray*}
First subtract the last row from the second last row,
and next factor out $(a_{2n-1}-a_{2n})$ from the second last row.
Then,
subtract $a_{2n}$ times the second last row from the last row.
Thus we obtain
\begin{eqnarray*}
&&V^{n}(X,\pmb{1};A,\pmb{1})\Big\rvert_{x_{2n}=x_{2n-1}}
=(a_{2n-1}-a_{2n})\\
&&\times\det\left(
\begin{cases}
a_{i}x_{i}^{n-j}
&\text{ if $1\leq i\leq2n-2$ and $1\leq j\leq n$,}\\
x_{i}^{2n-j}
&\text{ if $1\leq i\leq2n-2$ and $n+1\leq j\leq 2n$,}\\
x_{i}^{n-j}
&\text{ if $i=2n-1$ and $1\leq j\leq n$,}\\
0
&\text{ if $i=2n-1$ and $n+1\leq j\leq 2n$,}\\
0
&\text{ if $i=2n$ and $1\leq j\leq n$,}\\
x_{i}^{2n-j}
&\text{ if $i=2n$ and $n+1\leq j\leq 2n$.}
\end{cases}
\right)_{1\leq i,j\leq2n}.
\end{eqnarray*}
If we use our example,
then the matrix after this transformations looks like
\begin{eqnarray*}
V^{3}(X,\pmb{1};A,\pmb{1})\Big\rvert_{x_{6}=x_{5}}
=(a_5-a_6)
\begin{vmatrix}
a_{1}x_{1}^{2}&a_{1}x_{1}&a_{1}&
x_{1}^{2}&x_{1}&1\\
a_{2}x_{2}^{2}&a_{2}x_{2}&a_{2}&
x_{2}^{2}&x_{2}&1\\
a_{3}x_{3}^{2}&a_{3}x_{3}&a_{3}&
x_{3}^{2}&x_{3}&1\\
a_{4}x_{4}^{2}&a_{4}x_{4}&a_{4}&
x_{4}^{2}&x_{4}&1\\
x_{5}^{2}&x_{5}&1&
0&0&0\\
0&0&0&
x_{5}^{2}&x_{5}&1\\
\end{vmatrix}.
\end{eqnarray*}
Now subtract $x_{2n-1}$ times the second column from the first column,
then subtract $x_{2n-1}$ times the third column from the second column,
and so on,
until we subtract $x_{2n-1}$ times the $n$th column from the $(n-1)$th column.
Next subtract $x_{2n-1}$ times the $(n+2)$th column from the $(n+1)$th column,
then subtract $x_{2n-1}$ times the $(n+3)$th column from the $(n+2)$th column,
and so on,
until we subtract $x_{2n-1}$ times the $2n$th column from the $(2n-1)$th column.
If we expand the resulting determinant along the last two rows
and factor out $(x_i-x_{2n})$ from the $i$th row for $1\leq i\leq2n-2$,
then we obtain
\begin{eqnarray*}
&&V^{n}(X,\pmb{1};A,\pmb{1})\Big\rvert_{x_{2n}=x_{2n-1}}
=(-1)^{n+1}(a_{2n-1}-a_{2n})\prod_{i=1}^{2n-2}(x_i-x_{2n-1})\\
&&\qquad\qquad\times V^{n-1}(X^{(2n-1,2n)},\pmb{1}^{(2n-1,2n)};A^{(2n-1,2n)},\pmb{1}^{(2n-1,2n)}),
\end{eqnarray*}
and this proves our lemma.
$\Box$
\end{demo}
\begin{lemma}
\label{recurrence}
Let $n$ be a positive integer.
Let $X=(x_1,\dots,x_{2n})$, $A=(a_1,\dots,a_{2n})$ and $C=(c_1,\dots,c_{2n})$
be $2n$-tuples of variables.
Then the following identity holds.
\begin{align*}
&\sum_{k=1}^{2n-1}
\frac{\prod_{\scriptstyle i=1\atop\scriptstyle i\neq k}^{2n-1}(x_{k}-x_{i})}
{x_{k}-x_{2n}}
(a_{k}-a_{2n})
(c_{k}-c_{2n})
\\
&\times
V^{n-1}(X^{(k,2n)},\pmb{1}^{(k,2n)};A^{(k,2n)},\pmb{1}^{(k,2n)})
V^{n-1}(X^{(k,2n)},\pmb{1}^{(k,2n)};C^{(k,2n)},\pmb{1}^{(k,2n)})\\
&=\frac{V^{n}(X,\pmb{1};A,\pmb{1})V^{n}(X,\pmb{1};C,\pmb{1})}
{\prod_{i=1}^{2n-1}(x_{i}-x_{2n})}.
\end{align*}
Here $\pmb{1}$ denotes the $2n$-tuples $(1,\dots,1)$.
\end{lemma}
\begin{demo}{Proof}
It is enough to prove this lemma as an identity for a rational function in the complex variable $x_{2n}$.
Further we may assume $x_1$, $\dots$, $x_{2n-1}$ are distinct complex numbers.
Denote by $F(x_{2n})$ the left-hand side, and by $G(x_{2n})$ the right-hand side.
Under this assumption $F(x_{2n})$ and $G(x_{2n})$ have only simple poles as singularities,
and these simple poles reside at $x_{2n}=x_{k}$ for $1\leq k\leq 2n-1$.
First we want to show that
\begin{align*}
\Res_{x_{2n}=x_{k}}F(x_{2n})
=\Res_{x_{2n}=x_{k}}G(x_{2n}).
\end{align*}
The residue of the function $F(x_{2n})$ at $x_{2n}=x_{k}$ is
\begin{align*}
&\lim_{x_{2n}\rightarrow x_{k}}(x_{2n}-x_{k})F(x_{2n})\\
&=-(a_{k}-a_{2n})(c_{k}-c_{2n})
\prod_{\scriptstyle i=1\atop\scriptstyle i\neq k}^{2n-1}(x_{k}-x_{i})
\\
&\times
V^{n-1}(X^{(k,2n)},\pmb{1}^{(k,2n)};A^{(k,2n)},\pmb{1}^{(k,2n)})
V^{n-1}(X^{(k,2n)},\pmb{1}^{(k,2n)};C^{(k,2n)},\pmb{1}^{(k,2n)}).
\end{align*}
On the other hand,
the residue of the function $G(x_{2n})$ at $x_{2n}=x_{k}$ is
\begin{align*}
&\lim_{x_{2n}\rightarrow x_{k}}(x_{2n}-x_{k})G(x_{2n})
=-\frac{V^{n}(X,\pmb{1};A,\pmb{1})\Big\rvert_{x_{2n}=x_{k}}V^{n}(X,\pmb{1};C,\pmb{1})\Big\rvert_{x_{2n}=x_{k}}}
{\prod_{{i=1}\atop{i\neq k}}^{2n-1}(x_{i}-x_{k})}.
\end{align*}
By Lemma~\ref{substitution},
we see that $\Res_{x_{2n}=x_{k}}G(x_{2n})$ is equal to $\Res_{x_{2n}=x_{k}}F(x_{2n})$.
Thus it is shown that the principal part of $F(x_{2n})$ at each singularity coincides with that of $G(x_{2n})$.
Also it is clear that $\lim_{x_{2n}\rightarrow\infty}F(x_{2n})=\lim_{x_{2n}\rightarrow\infty}G(x_{2n})=0$.
Hence we have $F(x_{2n})=G(x_{2n})$.
$\Box$
\end{demo}
Now we are in the position to prove the first theorem in this section.
\begin{demo}{Proof of Theorem~\ref{fundamental}}
First we prove \thetag{\ref{eq:fundamental}}
when $Y=B=D=\pmb{1}=(1,\dots,1)$,
and deduce the general case to this special case.
Thus our first claim is
\begin{align}
\label{eq:special}
&\Pf\left[\frac{(a_{i}-a_{j})(c_{i}-c_{j})}{x_i-x_j}\right]_{1\leq i<j\leq 2n}
=\frac{V^{n}(X,\pmb{1};A,\pmb{1})V^{n}(X,\pmb{1};C,\pmb{1})}
{\prod_{1\leq i<j\leq2n}\left(x_{i}-x_{j}\right)}.
\end{align}
We proceed by induction on $n$.
When $n=1$, it is trivial since $V^{1}(X,\pmb{1};A,\pmb{1})=a_1-a_2$ and $V^{1}(X,\pmb{1};C,\pmb{1})=c_1-c_2$.
Assume $n\geq2$ and the identity holds up to $(n-1)$.
Expanding the Pfaffian along the last row/column and using the induction hypothesis,
we obtain
\begin{align*}
&\Pf\left[\frac{(a_{i}-a_{j})(c_{i}-c_{j})}{x_i-x_j}\right]_{1\leq i<j\leq 2n}\\
&=\sum_{k=1}^{2n-1}(-1)^{k-1}
\frac{(a_{k}-a_{2n})(c_{k}-c_{2n})}{x_{k}-x_{2n}}\\
&\times\frac{V^{n-1}(X^{(k,2n)},\pmb{1}^{(k,2n)};A^{(k,2n)},\pmb{1}^{(k,2n)})
V^{n-1}(X^{(k,2n)},\pmb{1}^{(k,2n)};C^{(k,2n)},\pmb{1}^{(k,2n)})}
{\displaystyle\prod_{{1\leq i<j\leq 2n-1}\atop{i,j\neq k}}(x_i-x_j)}.
\end{align*}
Substituting
$
\displaystyle\prod_{{1\leq i<j\leq 2n-1}\atop{i,j\neq k}}(x_i-x_j)
=(-1)^{k-1}
\frac{\displaystyle\prod_{1\leq i<j\leq 2n-1}(x_i-x_j)}
{\displaystyle\prod_{{i=1}\atop{i\neq k}}^{2n-1}(x_k-x_i)}
$
into this identity,
we have
\begin{align*}
&\Pf\left[\frac{(a_{i}-a_{j})(c_{i}-c_{j})}{x_i-x_j}\right]_{1\leq i<j\leq 2n}\\
&=\frac1{\displaystyle\prod_{1\leq i<j\leq 2n-1}(x_i-x_j)}
\sum_{k=1}^{2n-1}
(a_{k}-a_{2n})(c_{k}-c_{2n})
\cdot
\frac{\prod_{{i=1}\atop{i\neq k}}^{2n-1}(x_k-x_i)}
{x_{k}-x_{2n}}
\\&\times
V^{n-1}(X^{(k,2n)},\pmb{1}^{(k,2n)};A^{(k,2n)},\pmb{1}^{(k,2n)})
V^{n-1}(X^{(k,2n)},\pmb{1}^{(k,2n)};C^{(k,2n)},\pmb{1}^{(k,2n)}).
\end{align*}
Thus,
by Lemma~\ref{recurrence},
this equals
\begin{align*}
\frac{V^{n}(X,\pmb{1};A,\pmb{1})V^{n}(X,\pmb{1};C,\pmb{1})}
{\prod_{1\leq i<j\leq2n}(x_{i}-x_{j})},
\end{align*}
and this proves \thetag{\ref{eq:special}}.
The general identity \thetag{\ref{eq:fundamental}} is an easy consequence of \thetag{\ref{eq:special}}
by substituting $\frac{x_{i}}{y_{i}}$ into $x_{i}$,
$\frac{a_{i}}{b_{i}}$ into $b_{i}$ and $\frac{c_{i}}{d_{i}}$ into $c_{i}$
for $1\leq i\leq 2n$
in the both sides of \thetag{\ref{eq:special}}.
In fact,
if we write $X/Y$ for $(x_{1}/y_{1},\dots,x_{2n}/y_{2n})$,
$A/B$ for $(a_{1}/b_{1},\dots,a_{2n}/b_{2n})$,
and $C/D$ for $(c_{1}/d_{1},\dots,c_{2n}/d_{2n})$,
then we have
\begin{eqnarray*}
V^{n}(X/Y,\pmb{1};A/B,\pmb{1})
&=&\prod_{i=1}^{2n}b_{i}^{-1}\prod_{i=1}^{2n}y_{i}^{-n+1}
\cdot V^{n}(X,Y;A,B)\\
V^{n}(X/Y,\pmb{1};C/D,\pmb{1})
&=&\prod_{i=1}^{2n}d_{i}^{-1}\prod_{i=1}^{2n}y_{i}^{-n+1}
\cdot V^{n}(X,Y;C,D).
\end{eqnarray*}
Since the denominator of the right-hand side of \thetag{\ref{eq:special}} is
\begin{eqnarray*}
\prod_{1\leq i<j\leq2n}(x_{i}/y_{i}-x_{j}/y_{j})
=\prod_{i=1}^{2n}y_{i}^{-2n+1}
\prod_{1\leq i<j\leq2n}
\begin{vmatrix}
x_i&y_i\\
x_j&y_j
\end{vmatrix},
\end{eqnarray*}
the right-hand side of \thetag{\ref{eq:special}} becomes
\begin{align*}
&\frac{V^{n}(X/Y,\pmb{1};A/B,\pmb{1})V^{n}(X/Y,\pmb{1};C/D,\pmb{1})}
{\prod_{1\leq i<j\leq2n}(x_{i}/y_{i}-x_{j}/y_{j})}\\
&=
\prod_{i=1}^{2n}b_{i}^{-1}\prod_{i=1}^{2n}d_{i}^{-1}\prod_{i=1}^{2n}y_{i}\cdot
\frac{V^{n}(X,Y;A,B)V^{n}(X,Y;C,D)}
{\prod_{1\leq i<j\leq2n}\begin{vmatrix}x_{i}&y_{i}\\x_{j}&y_{j}\end{vmatrix}}
\end{align*}
On the other hand,
the left-hand side of \thetag{\ref{eq:special}} becomes
\begin{align*}
\Pf\left[
\frac{(\frac{a_{i}}{b_{i}}-\frac{a_{j}}{b_{j}})(\frac{c_{i}}{d_{i}}-\frac{c_{j}}{d_{j}})}
{\frac{x_{i}}{y_{i}}-\frac{x_{j}}{y_{j}}}
\right]_{1\leq i<j\leq2n}
=\prod_{i=1}^{2n}\left(b_{i}^{-1}d_{i}^{-1}y_{i}\right)\cdot
\Pf\left[
\frac{\begin{vmatrix}a_{i}&b_{i}\\a_{j}&b_{j}\end{vmatrix}
\cdot\begin{vmatrix}c_{i}&d_{i}\\c_{j}&d_{j}\end{vmatrix}}
{\begin{vmatrix}x_{i}&y_{i}\\x_{j}&y_{j}\end{vmatrix}}
\right]_{1\leq i<j\leq2n}.
\end{align*}
and this proves \thetag{\ref{eq:fundamental}}.
$\Box$
\end{demo}
Next we give a proof of Proposition~\ref{prop:subs}.
\begin{demo}{Proof of Proposition~\ref{prop:subs}}
From definition,
$V^{n}(X,\pmb{1}+tX^2;X,\pmb{1})$ is equal to
\begin{eqnarray*}
\det\left(
\begin{cases}
x_{i}^{n-j+1}(1+tx_{i}^{2})^{j-1}
&\text{ if $1\leq j\leq n$,}\\
x_{i}^{2n-j}(1+tx_{i}^{2})^{j-n-1}
&\text{ if $n+1\leq j\leq 2n$.}
\end{cases}
\right)_{1\leq i,j\leq2n}.
\end{eqnarray*}
For example,
if $n=3$,
then $V^{3}(X,\pmb{1}+tX^2;X,\pmb{1})$
is equal to
\begin{eqnarray*}
\begin{vmatrix}
x_{1}^{3}&x_{1}^{2}(1+tx_{1}^{2})&x_{1}(1+tx_{1}^{2})^2&
x_{1}^{2}&    x_{1}(1+tx_{1}^{2})&(1+tx_{1}^{2})^2\\
x_{2}^{3}&x_{2}^{2}(1+tx_{2}^{2})&x_{2}(1+tx_{2}^{2})^2&
x_{2}^{2}&    x_{2}(1+tx_{2}^{2})&(1+tx_{2}^{2})^2\\
x_{3}^{3}&x_{3}^{2}(1+tx_{3}^{2})&x_{3}(1+tx_{3}^{2})^2&
x_{3}^{2}&    x_{3}(1+tx_{3}^{2})&(1+tx_{3}^{2})^2\\
x_{4}^{3}&x_{4}^{2}(1+tx_{4}^{2})&x_{4}(1+tx_{4}^{2})^2&
x_{4}^{2}&    x_{4}(1+tx_{4}^{2})&(1+tx_{4}^{2})^2\\
x_{5}^{3}&x_{5}^{2}(1+tx_{5}^{2})&x_{5}(1+tx_{5}^{2})^2&
x_{5}^{2}&    x_{5}(1+tx_{5}^{2})&(1+tx_{5}^{2})^2\\
x_{6}^{3}&x_{6}^{2}(1+tx_{6}^{2})&x_{6}(1+tx_{6}^{2})^2&
x_{6}^{2}&    x_{6}(1+tx_{6}^{2})&(1+tx_{6}^{2})^2\\
\end{vmatrix}.
\end{eqnarray*}
By expanding $(1+tx_{i}^{2})^{j-1}$ and $(1+tx_{i}^{2})^{j-n-1}$,
and performing appropriate elementary column transformations,
this determinant becomes
\begin{eqnarray*}
\det\left(
\begin{cases}
x_{i}^{n}
&\text{ if $j=1$,}\\
x_{i}^{n+1-j}+t^{j-1}x_{i}^{j+n-1}
&\text{ if $2\leq j\leq n$,}\\
x_{i}^{n-1}
&\text{ if $j=2n+1$,}\\
x_{i}^{2n-j}+t^{j-n-1}x_{i}^{j-2}
&\text{ if $n+2\leq j\leq 2n$.}
\end{cases}
\right)_{1\leq i,j\leq2n}.
\end{eqnarray*}
For example,
if $n=3$,
then this determinant equals
\begin{eqnarray*}
\begin{vmatrix}
x_{1}^{3}&x_{1}^{2}+tx_{1}^{4}&x_{1}+t^2x_{1}^{5}&
x_{1}^{2}&    x_{1}+tx_{1}^{3}&    1+t^2x_{1}^{4}\\
x_{2}^{3}&x_{2}^{2}+tx_{2}^{4}&x_{2}+t^2x_{2}^{5}&
x_{2}^{2}&    x_{2}+tx_{2}^{3}&    1+t^2x_{2}^{4}\\
x_{3}^{3}&x_{3}^{2}+tx_{3}^{4}&x_{3}+t^2x_{3}^{5}&
x_{3}^{2}&    x_{3}+tx_{3}^{3}&    1+t^2x_{3}^{4}\\
x_{4}^{3}&x_{4}^{2}+tx_{4}^{4}&x_{4}+t^2x_{4}^{5}&
x_{4}^{2}&    x_{4}+tx_{4}^{3}&    1+t^2x_{4}^{4}\\
x_{5}^{3}&x_{5}^{2}+tx_{5}^{4}&x_{5}+t^2x_{5}^{5}&
x_{5}^{2}&    x_{5}+tx_{5}^{3}&    1+t^2x_{5}^{4}\\
x_{6}^{3}&x_{6}^{2}+tx_{6}^{4}&x_{6}+t^2x_{6}^{5}&
x_{6}^{2}&    x_{6}+tx_{6}^{3}&    1+t^2x_{6}^{4}
\end{vmatrix}.
\end{eqnarray*}
We subtract $t$ times the first column from the $(n+2)$th column,
and subtract the $(n+1)$th column from the second column.
Then we subtract $t$ times the second column from the $(n+3)$th column,
and subtract the $(n+2)$th column from the third column.
We continue this transformation
until we subtract $t$ times the $(n-1)$th column from the $2n$th column,
and subtract the $(2n-1)$th column from the $n$th column.
Thus we obtain
\begin{eqnarray*}
V^{n}(X,\pmb{1}+tX^2;X,\pmb{1})
=
\left(\begin{cases}
t^{j-1}x_{i}^{j+n-1}
&\text{ if $1\leq j\leq n$,}\\
x_{i}^{2n-j}
&\text{ if $n+1\leq j\leq 2n$,}\\
\end{cases}
\right)_{1\leq i,j\leq2n}.
\end{eqnarray*}
If we illustrate by the above example,
then this determinant looks like
\begin{eqnarray*}
\begin{vmatrix}
x_{1}^{3}&tx_{1}^{4}&t^2x_{1}^{5}&
x_{1}^{2}&     x_{1}&    1\\
x_{2}^{3}&tx_{2}^{4}&t^2x_{2}^{5}&
x_{2}^{2}&     x_{2}&    1\\
x_{3}^{3}&tx_{3}^{4}&t^2x_{3}^{5}&
x_{3}^{2}&     x_{3}&    1\\
x_{4}^{3}&tx_{4}^{4}&t^2x_{4}^{5}&
x_{4}^{2}&     x_{4}&    1\\
x_{5}^{3}&tx_{5}^{4}&t^2x_{5}^{5}&
x_{5}^{2}&     x_{5}&    1\\
x_{6}^{3}&tx_{6}^{4}&t^2x_{6}^{5}&
x_{6}^{2}&     x_{6}&    1
\end{vmatrix}.
\end{eqnarray*}
By the Vandermonde determinant,
we can easily conclude that
\begin{eqnarray*}
V^{n}(X,\pmb{1}+tX^2;X,\pmb{1})
=(-1)^{\binom{n}2}t^{0+1+2+\dots+(n-1)}\prod_{1\leq i<j\leq2n} (x_i-x_j).
\end{eqnarray*}
This completes the proof.
$\Box$
\end{demo}
Corollary~\ref{cor:pfaff-det} is an immediate consequence of Theorem~\ref{fundamental} and Proposition~\ref{prop:subs}.
So, in the rest of this section we derive Corollary~\ref{cor:Sunduist} from Theorem~\ref{fundamental}.
The following proposition is an immediate consequence of the Laplace expansion formula
and the Littlewood formula:
\begin{equation}
\label{Littlewood}
\sum_{\lambda\in\mathcal{P}_n}s_{\lambda}(x_1,\dots,x_n)
=\prod_{1\leq i<j\leq n}(1+x_ix_j),
\end{equation}
(see \cite{Ma} I.5. Ex.9).
Let $\binom{S}{r}$ denote the set of all $r$-element subsets of $S$
for any finite set $S$ and a non-negative integer $r$.
\begin{prop}
\label{Laplace}
Let $n$ be a non-negative integer.
Let $X=(x_1,\dots,x_{2n})$ and $A=(a_1,\dots,a_{2n})$
be $2n$-tuples of variables.
Then we have
\begin{align*}
U^{n}(X;A)
&=\sum_{I\in\binom{[2n]}{n}}(-1)^{|I|+\binom{n+1}{2}}a_{I}
\prod_{{i,j\in I}\atop{i<j}}(x_{i}-x_{j})(1+x_{i}x_{j})\\
&\qquad\qquad\times
\prod_{{i,j\in{I^c}}\atop{i<j}}(x_{i}-x_{j})(1+x_{i}x_{j}),
\end{align*}
where
$I^c=[2n]\setminus I$ is the complementary set of $I$ in $[2n]$
and $a_{I}=\prod_{i\in I}a_i$ for any subset $I\subseteq[2n]$.
\end{prop}
\begin{demo}{Proof}
We write $\Delta(X_{I})=\prod_{{i,j\in I}\atop{i<j}}(x_i-x_j)$.
By the Laplace expansion formula
we have
\begin{align*}
a_{\lambda,\mu}(X;A)
=\sum_{I\in\binom{[2n]}{n}}(-1)^{|I|+\binom{n+1}2}
a_I\Delta(X_{I})\Delta(X_{I^c})s_{\lambda}(X_{I})s_{\mu}(X_{I^c}),
\end{align*}
where $X_{I}$ stands for the $n$-tuple of variables with index set $I$ and $X_{I^c}$ stands for the $n$-tuple of variables with index set $I^c$,
and $|I|=\sum_{i\in I}i$.
By \thetag{\ref{Littlewood}}
we easily obtain the desired formula.
$\Box$
\end{demo}
\begin{demo}{Proof of Theorem~\ref{det_exp}}
As before
we apply the Laplace expansion formula to $V^{n}(X,\pmb{1}-X^2;A,\pmb{1})$ to obtain
\begin{eqnarray*}
&&V^{n}(X,\pmb{1}-X^2;A,\pmb{1})
=\sum_{I\in\binom{[2n]}{n}}
(-1)^{|I|+\binom{n+1}2}a_I
\det\Delta^{I}(M)
\det\Delta^{I^c}(M),
\end{eqnarray*}
where $M=\left(x_{i}^{n-j}(1-x_{i}^2)^{j-1}\right)_{1\leq i\leq 2n,\,1\leq j\leq n}$.
By the Vandermonde determinant,
we obtain
\begin{align*}
\det\Delta^{I}(M)
&=\prod_{{i,j\in I}\atop{i<j}}
\left\{x_{i}(1-x_{j}^{2})-x_{j}(1-x_{i}^{2})\right\}\\
&=\prod_{{i,j\in I}\atop{i<j}}(x_i-x_j)(1+x_ix_j).
\end{align*}
Thus we conclude that
\begin{eqnarray*}
V^{n}(X,\pmb{1}-X^2;A,\pmb{1})
&=&\sum_{I\in\binom{[2n]}{n}}
(-1)^{|I|+\binom{n+1}2}a_I
\prod_{{i,j\in I}\atop{i<j}}(x_i-x_j)(1+x_ix_j)\\
&&\qquad\times\prod_{{i,j\in I^c}\atop{i<j}}(x_i-x_j)(1+x_ix_j).
\end{eqnarray*}
Thus,
by Proposition~\ref{Laplace},
we complete our proof.
$\Box$
\end{demo}
\begin{demo}{Proof of Corollary~\ref{cor:Sunduist}}

If we substitute $x_i^2$ into $x_i$, $x_i$ into $c_i$,
and $1$ into $y_i$, $b_i$ and $d_i$ for $1\leq i\leq 2n$ in \thetag{\ref{eq:fundamental}},
then we obtain
\begin{align*}
&\Pf\left[\frac{a_i-a_j}{x_i+x_j}\right]_{1\leq i<j\leq 2n}
=\frac{V^{n}(X^2,\pmb{1};A,\pmb{1})V^{n}(X^2,\pmb{1};X,\pmb{1})}{\prod_{1\leq i<j\leq 2n}(x_i^2-x_j^2)},
\end{align*}
where $X^2=(x_{1}^2,\dots,x_{2n}^2)$.
By Proposition~\ref{prop:subs}
\[
V^{n}(X^2,\pmb{1};X,\pmb{1})
=(-1)^{\binom{n}2}\prod_{1\leq i<j\leq2n}(x_i-x_j)
\]
and
this proves \thetag{\ref{Sundquist1}}.

Next we state a proof of another identity.
Substituting $t=-1$ and $b_i=1$ for $1\leq i\leq 2n$ in \thetag{\ref{eq:cauchy}},
we obtain
\begin{align*}
&\Pf\left[\frac{a_i-a_j}
{1+x_ix_j}\right]_{1\leq i<j\leq 2n}
=\frac{V^{n}(X,\pmb{1}-X^2;A,\pmb{1})}{\prod_{1\leq i<j\leq 2n}(1+x_ix_j)}.
\end{align*}
Thus, by Theorem~\ref{det_exp},
we immediately obtain the desired identity.
This completes the proof.
$\Box$
\end{demo}

%
%
\section{A Proof of Stanley's Open Problem}

The key idea of our proof is the following proposition,
which the reader can find in \cite{Sta1}, Exercise~7.7, or \cite{Ste2}, Section~3.
\begin{prop}
\label{key_lemma}
Let $f(x_1,x_2,\dots)$ be a symmetric function with infinite variables.
Then $f\in{\Bbb{Q}}[p_\lambda:\text{ all parts $\lambda_i>0$ are odd}]$
if and only if
\begin{eqnarray*}
f(t,-t,x_1,x_2,\dots)=f(x_1,x_2,\dots).
\Box
\end{eqnarray*}
\end{prop}
\noindent
Our strategy is simple.
If we set $v_n(X_{2n})$ to be
\begin{align}
\label{eq:v_n}
\log z_n(X_{2n})
-\sum_{k\geq1}\frac{1}{2k}a^k(b^k-c^k)p_{2k}(X_{2n})
-\sum_{k\geq1}\frac{1}{4k}a^kb^kc^kd^kp_{2k}(X_{2n})^2
\end{align}
then we claim it satisfies
\begin{align}
\label{eq:result}
v_{n+1}(t,-t,X_{2n})=v_n(X_{2n}).
\end{align}
This will eventually prove Theorem~\ref{conj:general}.
As an immediate consequence of \thetag{\ref{eq:general}}, \thetag{\ref{eq:entries}} and \thetag{\ref{eq:cauchy}},
we obtain the following theorem:
\begin{theorem}
\label{cor:det}
Let $X=(x_{1},\dots,x_{2n})$ be a $2n$-tuple of variables.
Then
\begin{align}
\label{eq:det_zn}
z_{n}(X_{2n})=(-1)^{\binom{n}2}\frac
{V^{n}(X^2,\pmb1+abcd X^4;X+aX^2,\pmb1-a(b+c)X^2-abcX^3)}
{\prod_{i=1}^{2n}(1-abx_{i}^2)\prod_{1\leq i<j\leq2n}(x_{i}-x_{j})(1-abcdx_{i}^2x_{j}^2)},
\end{align}
where $X^2=(x_{1}^2,\dots,x_{2n}^2)$, $\pmb1+abcd X^4=(1+abcdx_{1}^4,\dots,1+abcdx_{2n}^4)$, $X+aX^2=(x_1+ax_1^2,\dots,x_{2n}+ax_{2n}^2)$ and $\pmb1-a(b+c)X^2-abcX^3=(1-a(b+c)x_{1}^2-abcx_{1}^3,\dots,1-a(b+c)x_{2n}^2-abcx_{2n}^3)$.
$\Box$
\end{theorem}
The \thetag{\ref{eq:det_zn}} is key expression to prove that $v_n(X_{2n})$ satisfies \thetag{\ref{eq:result}}.
Once one knows \thetag{\ref{eq:det_zn}},
then it is straight forward computation to prove Stanley's open problem.
The following proposition is the first step.
\begin{prop}
\label{prop:numer}
Let $X=(x_{1},\dots,x_{2n})$ be a $2n$-tuple of variables.
Put
\[
f_{n}(X_{2n})=V^{n}(X^2,\pmb1+abcd X^4;X+aX^2,\pmb1-a(b+c)X^2-abcX^3).
\]
Then
$f_{n}(X_{2n})$ satisfies
\begin{align}
\label{eq:f_n}
&f_{n+1}(t,-t,X_{2n})
\nonumber\\
&=(-1)^{n}
2t(1-abt^2)(1-act^2)
\prod_{i=1}^{2n}(t^2-x_{i}^2)
\prod_{i=1}^{2n}(1-abcdt^2x_{i}^2)
\cdot f_{n}(X_{2n}).
\end{align}
\end{prop}
\begin{demo}{Proof}
First,
we put $\xi_{i}=x_{i}^2$, $\eta_{i}=1+abcdx_{i}^4$,
$\alpha_{i}=x_{i}+ax_{i}^2$, $\beta_{i}=1-a(b+c)x_{i}^2-abcx_{i}^3$
and $\zeta_{i}=\xi_{i}^{-1}\eta_{i}=x_{i}^{-2}+abcdx_{i}^2$ for $1\leq i\leq2n$.
Then
\begin{align*}
f_{n+1}(X_{2n+2})
&=\det\begin{pmatrix}
\begin{cases}
\alpha_{i}\xi_{i}^{n+1-j}\eta_{i}^{j-1}
&\text{ if $1\leq j\leq n+1$,}\\
\beta_{i}\xi_{i}^{2n+2-j}\eta_{i}^{j-n-2}
&\text{ if $n+2\leq j\leq 2n+2$.}
\end{cases}
\end{pmatrix}_{1\leq i,j\leq2n+2},\\
&=\prod_{i=1}^{2n}\xi_{i}^{n}\cdot
\det\begin{pmatrix}
\begin{cases}
\alpha_{i}\zeta_{i}^{j-1}
&\text{ if $1\leq j\leq n+1$,}\\
\beta_{i}\zeta_{i}^{j-n-2}
&\text{ if $n+2\leq j\leq 2n+2$.}
\end{cases}
\end{pmatrix}_{1\leq i,j\leq2n+2}.
\end{align*}
For example,
if $n=2$
then $f_{3}(X_{6})$ looks as follows:
\begin{equation*}
\prod_{i=1}^{2n}\xi_{i}^{2}\cdot
\begin{vmatrix}
\alpha_{1}&\alpha_{1}\zeta_{1}&\alpha_{1}\zeta_{1}^2&\beta_{1}&\beta_{1}\zeta_{1}&\beta_{1}\zeta_{1}^2\\
\alpha_{2}&\alpha_{2}\zeta_{2}&\alpha_{2}\zeta_{2}^2&\beta_{2}&\beta_{2}\zeta_{2}&\beta_{2}\zeta_{2}^2\\
\alpha_{3}&\alpha_{3}\zeta_{3}&\alpha_{3}\zeta_{3}^2&\beta_{3}&\beta_{3}\zeta_{3}&\beta_{3}\zeta_{3}^2\\
\alpha_{4}&\alpha_{4}\zeta_{4}&\alpha_{4}\zeta_{4}^2&\beta_{4}&\beta_{4}\zeta_{4}&\beta_{4}\zeta_{4}^2\\
\alpha_{5}&\alpha_{5}\zeta_{5}&\alpha_{5}\zeta_{5}^2&\beta_{5}&\beta_{5}\zeta_{5}&\beta_{5}\zeta_{5}^2\\
\alpha_{6}&\alpha_{6}\zeta_{6}&\alpha_{6}\zeta_{6}^2&\beta_{6}&\beta_{6}\zeta_{6}&\beta_{6}\zeta_{6}^2\\
\end{vmatrix}.
\end{equation*}
Now we subtract $\zeta_{1}$ times the $n$th column from the $(n+1)$th column,
then subtract $\zeta_{1}$ times the $(n-1)$th column from the $n$th column,
and so on,
until we subtract $\zeta_{1}$ times the first column from the second column.
Next 
we subtract $\zeta_{1}$ times the $(2n+1)$th column from the $(2n+2)$th column,
then subtract $\zeta_{1}$ times the $2n$th column from the $(2n+1)$th column,
and so on,
until we subtract $\zeta_{1}$ times the $(n+2)$th column from the $(n+3)$th column.
Thus we obtain $f_{n+1}(X_{2n+2})$ is equal to
\begin{equation*}
\prod_{i=1}^{2n}\xi_{i}^{n}\cdot
\det\left(\begin{cases}
\alpha_{1}
&\text{ if $i=1$ and $j=1$,}\\
\beta_{1}
&\text{ if $i=1$ and $j=n+2$,}\\
0
&\text{ if $i=1$ and $j\neq1,n+2$,}\\
\alpha_{i}\zeta_{i}^{j-2}(\zeta_{i}-\zeta_{1})
&\text{ if $i\geq2$ and $1\leq j\leq n+1$,}\\
\beta_{i}\zeta_{i}^{j-n-3}(\zeta_{i}-\zeta_{1})
&\text{ if $i\geq2$ and $n+2\leq j\leq 2n+2$.}
\end{cases}
\right)_{1\leq i,j\leq2n+2}.
\end{equation*}
If we illustrate in the above example,
then this determinant looks
\begin{equation*}
\prod_{i=1}^{2n}\xi_{i}^{2}\cdot
\begin{vmatrix}
\alpha_{1}&0&0&\beta_{1}&0&0\\
\alpha_{2}&\alpha_{2}(\zeta_{2}-\zeta_{1})&\alpha_{2}\zeta_{2}(\zeta_{2}-\zeta_{1})&\beta_{2}&\beta_{2}(\zeta_{2}-\zeta_{1})&\beta_{2}\zeta_{2}(\zeta_{2}-\zeta_{1})\\
\alpha_{3}&\alpha_{3}(\zeta_{3}-\zeta_{1})&\alpha_{3}\zeta_{3}(\zeta_{3}-\zeta_{1})&\beta_{3}&\beta_{3}(\zeta_{3}-\zeta_{1})&\beta_{3}\zeta_{3}(\zeta_{3}-\zeta_{1})\\
\alpha_{4}&\alpha_{4}(\zeta_{4}-\zeta_{1})&\alpha_{4}\zeta_{4}(\zeta_{4}-\zeta_{1})&\beta_{4}&\beta_{4}(\zeta_{4}-\zeta_{1})&\beta_{4}\zeta_{4}(\zeta_{4}-\zeta_{1})\\
\alpha_{5}&\alpha_{5}(\zeta_{5}-\zeta_{1})&\alpha_{5}\zeta_{5}(\zeta_{5}-\zeta_{1})&\beta_{5}&\beta_{5}(\zeta_{5}-\zeta_{1})&\beta_{5}\zeta_{5}(\zeta_{5}-\zeta_{1})\\
\alpha_{6}&\alpha_{6}(\zeta_{6}-\zeta_{1})&\alpha_{6}\zeta_{6}(\zeta_{6}-\zeta_{1})&\beta_{6}&\beta_{6}(\zeta_{6}-\zeta_{1})&\beta_{6}\zeta_{6}(\zeta_{6}-\zeta_{1})
\end{vmatrix}.
\end{equation*}
Here,
if we assume $\xi_{1}=\xi_{2}$ and $\zeta_{1}=\zeta_{2}$ hold,
then $f_{n+1}(X_{2n+2})$ is equal to
\begin{eqnarray*}
(-1)^{n}
\begin{vmatrix}
\alpha_{1}&\beta_{1}\\
\alpha_{2}&\beta_{2}
\end{vmatrix}
\cdot
\prod_{i=3}^{2n+2}
\begin{vmatrix}
\xi_{1}&\eta_{1}\\
\xi_{i}&\eta_{i}
\end{vmatrix}
\cdot
f_{n}(x_{3},\dots,x_{2n+2}).
\end{eqnarray*}
Now we substitute $X_{2n+2}=(t,-t,X_{2n})$ into this identity,
then,
since we have $\xi_{1}=\xi_{2}=t^{2}$, $\zeta_{1}=\zeta_{2}=t^{-2}+abcdt^2$
and
\begin{align*}
&\begin{vmatrix}
\xi_{1}&\eta_{1}\\
\xi_{i}&\eta_{i}
\end{vmatrix}
=(t^2-x_{i}^2)(1-abcdx_{i}^2),
\\
&\begin{vmatrix}
\alpha_{1}&\beta_{1}\\
\alpha_{2}&\beta_{2}
\end{vmatrix}
=2t(1-abt^2)(1-act^2),
\end{align*}
thus we obtain
\begin{align*}
f_{n+1}(t,-t,X_{2n})
&=(-1)^{n}
\cdot 2t(1-abt^2)(1-act^2)\\
&\quad\times\prod_{i=1}^{2n}(t^2-x_{i}^2)(1-abcdt^2x_{i}^2)
\cdot f_{n}(X_{2n}).
\end{align*}
This proves our proposition.
$\Box$
\end{demo}
\begin{prop}
Let $X=(x_{1},\dots,x_{2n})$ be a $2n$-tuple of variables.
Then
\begin{align}
\label{eq:z_n}
z_{n+1}(t,-t,X_{2n})
=\frac{1-act^2}{(1-abt^2)(1-abcdt^4)\prod_{i=1}^{2n}(1-abcdt^2x_{i}^2)}
z_{n}(X_{2n}).
\end{align}
\end{prop}
\begin{demo}{Proof}
By Theorem~\ref{cor:det}
we have
\begin{eqnarray*}
z_{n}(X_{2n})
=(-1)^{\binom{n}2}\frac{f_{n}(X_{2n})}
{\prod_{i=1}^{2n}(1-abx_{i}^2)\prod_{1\leq i<j\leq2n}(x_{i}-x_{j})(1-abcdx_{i}^2x_{j}^2)}.
\end{eqnarray*}
This implies
\begin{eqnarray*}
z_{n+1}(t,-t,X_{2n})
&=(-1)^{\binom{n+1}2}\frac{1}{2t(1-abt^2)^2(1-abcdt^4)\cdot\prod_{i=1}^{2n}(t^2-x_{i}^2)(1-abcdt^2x_{i}^2)^2}
\\
&\times\frac{f_{n+1}(t,-t,X_{2n})}
{\prod_{i=1}^{2n}(1-abx_{i}^2)\prod_{1\leq i<j\leq2n}(x_{i}-x_{j})(1-abcdx_{i}^2x_{j}^2)}.
\end{eqnarray*}
Thus, substituting \thetag{\ref{eq:f_n}},
we obtain the desired identity.
$\Box$
\end{demo}
Now we are in the position to complete our proof of Stanley's open problem.
\begin{demo}{Proof of Theorem~\ref{conj:general}}
\thetag{\ref{eq:z_n}} immediately implies
\begin{eqnarray}
\log{z_{n+1}(t,-t,X_{2n})}
&=\log{z_{n}(X_{2n})}+\log\frac1{1-abt^2}-\log\frac1{1-act^2}
\nonumber\\
&+\log\frac1{1-abcdt^4}+\sum_{i=1}^{2n}\log\frac1{1-abcdt^2x_i^2}.
\label{eq:logz}
\end{eqnarray}
On the other hand,
$p_{2k}(t,-t,X_{2n})=2t^{2k}+\sum_{i=1}^{2n}x_{i}^{2k}$
implies
\begin{align*}
&\sum_{k\geq1}\frac{a^n(b^n-c^n)}{2k}p_{2k}(t,-t,X_{2n})
\\
&=\sum_{k\geq1}\frac{a^n(b^n-c^n)}{2k}p_{2k}(X_{2n})
+\log\frac1{1-abt^2}-\log\frac1{1-act^2},
\\
&\sum_{k\geq1}\frac{a^nb^nc^nd^n}{4k}p_{2k}(t,-t,X_{2n})^2
\\
&=\sum_{k\geq1}\frac{a^nb^nc^nd^n}{4k}p_{2k}(X_{2n})^2+\log\frac1{1-abcdt^4}+\sum_{i=1}^{2n}\log\frac1{1-abcdt^2x_i^2}.
\end{align*}
Thus, putting $v_{n}(X_{2n})$ as in \thetag{\ref{eq:v_n}},
we easily find $v_{n}(X_{2n})$ satisfies \thetag{\ref{eq:result}} from \thetag{\ref{eq:logz}}.
This completes our proof of Theorem~\ref{conj:general}.
$\Box$
\end{demo}

%
%
\section{Corollaries}

The author tried to find an analogous formula when the sum runs over all distinct partitions
by computer experiments using Stembridge's SF package
(cf. \cite{B1} and \cite{B2}).
But the author could not find any conceivable formula when the sum runs over all distinct partitions,
and, instead, found the following formula involving the big Schur functions
and certain symmetric functions arising from the Macdonald polynomials as byproducts.
But,
later,
Prof. R.~Stanley and Prof. A.~Lascoux independently pointed out 
these are derived from Theorem~\ref{conj:general} as corollaries.
Let $S_{\lambda}(x;t)=\det\left(q_{\lambda_i-i+j}(x;t)\right)_{1\leq i,j\leq\ell(\lambda)}$ denote the big Schur function corresponding to the partition $\lambda$,
where $q_r(x;t)=Q_{(r)}(x;t)$ denotes the Hall-Littlewood function
(See \cite{Ma}, III, sec.2).

\begin{corollary}
\label{conj:BigSchur}
Let
\begin{equation*}
Z(x;t)=\sum_{\lambda}\omega(\lambda)S_{\lambda}(x;t),
\end{equation*}
Here the sum runs over all partitions $\lambda$.
Then we have
\begin{align}
\log Z(x;t)-\sum_{n\geq1}\frac1{2n}a^n(b^n-c^n)(1-t^{2n})p_{2n}
-\sum_{n\geq1}\frac1{4n}a^nb^nc^nd^n(1-t^{2n})^2p_{2n}^2\nonumber\\
\in\Bbb{Q}[[p_1,p_3,p_5,\dots]].
\label{eq:BigSchur}
\end{align}
\end{corollary}
\begin{demo}{Proof}
This proof was originally suggested by Prof. R.~Stanley.
Let $\Lambda_x$ denote the ring of symmetric functions in countably many variables $x_1$, $x_2$, $\dots$.
(For details see \cite{Ma}, I, sec.2).
Let $\theta_x$ be the ring homomorphism $\Lambda_x \to \Lambda_x[t]$
taking $h_n(x)$ to $q_n(x;t)$. 
By the Jacobi-Trudi identity we have 
\begin{equation}
\label{eq:JT}
s_{\lambda}(x)=\det(h_{\lambda_i-i+j}(x)).
\end{equation}
(cf. \cite{Ma}, I, sec.2 (3.4).)
Applying $\theta_x$ to the both sides,
and using the definition of the big Schur $S_{\lambda}(x;t)=\det(q_{\lambda_i-i+j}(x;t))$
(\cite{Ma}, III, sec.4, (4.5)),
we obtain
\[
\theta_x(s_\lambda(x)) = S_\lambda(x;t).
\] 
By taking logarithms of
\[
   \sum_{\lambda} S_\lambda(x;t)s_\lambda(y) =
                 \prod_{i,j\geq1}\frac{1-tx_iy_j}{1-x_iy_j},
\]
(\cite{Ma}, III, sec.4, (4.7)),
 the product on the right-hand side is
\[
   \exp \sum_{n\geq 1} \frac1{n}(1-t^n)p_n(x)p_n(y).
\]
Similarly,
by taking logarithms of the right-hand side of
\begin{equation}
\label{Schur-prod}
\sum_{\lambda} s_\lambda(x)s_\lambda(y) =
                 \prod_{i,j\geq1}(1-x_iy_j)^{-1}
\end{equation}
(\cite{Ma}, I, sec4, (4.3)),
we have
\begin{equation*}
\label{eq:Schur-prod}
      \sum_\lambda s_\lambda(x)s_\lambda(y)
      = \exp \sum_{n\geq 1} \frac1{n}p_n(x)p_n(y), 
\end{equation*}
and it follows that $\theta(p_n) = (1-t^n)p_n$. 
The identity \thetag{\ref{eq:BigSchur}} now follows
by applying $\theta$ to equation \thetag{\ref{conj:general}}. 
$\Box$
\end{demo}
This corollary is generalized to the two parameter polynomials defined by I.~G.~Macdonald.
Define
\begin{eqnarray*}
T_{\lambda}(x;q,t)=\det\left(Q_{(\lambda_i-i+j)}(x;q,t)\right)_{1\leq i,j\leq\ell(\lambda)}
\end{eqnarray*}
where $Q_{\lambda}(x;q,t)$ stands for the Macdonald polynomial corresponding to the partition $\lambda$,
and $Q_{(r)}(x;q,t)$ is the one corresponding to the one row partition $(r)$
(See \cite{Ma}, IV, sec.4).
\begin{corollary}
\label{conj:Macdonald}
Let
\begin{equation*}
Z(x;q,t)=\sum_{\lambda}\omega(\lambda)T_{\lambda}(x;q,t),
\end{equation*}
Here the sum runs over all partitions $\lambda$.
Then we have
\begin{align}
\log Z(x;q,t)-\sum_{n\geq1}\frac1{2n}a^n(b^n-c^n)\frac{1-t^{2n}}{1-q^{2n}}p_{2n}
-\sum_{n\geq1}\frac1{4n}a^nb^nc^nd^n\frac{(1-t^{2n})^2}{(1-q^{2n})^2}p_{2n}^2
\nonumber\\
\in\Bbb{Q}[[p_1,p_3,p_5,\dots]].
\end{align}
\end{corollary}
\begin{demo}{Proof}
The proof proceeds almost parallel to that of Corollary~\ref{conj:BigSchur}
except that we define the ring homomorphism $\theta:\Lambda \to \Lambda(t,q)$
by $\theta(h_n)=g_n(x;q,t)$.
Here we write $g_n(x;q,t)=Q_{(n)}(x;q,t)$,
following the notation in \cite{Ma}.
Since
\[
\sum_{n\geq0}g_n(x;q,t)y^n
=\prod_{i\geq1}
\frac{(tx_iy;q)_{\infty}}{(x_iy;q)_{\infty}},
\]
(\cite{Ma}, VI, sec.2 (2.8)),
if we introduce a set of fictitious variables $\xi_i$ by
\[
\frac{(tx_iy;q)_{\infty}}{(x_iy;q)_{\infty}}=\prod_{i\geq1}(1-\xi_iy)^{-1},
\]
then we have $g_{r}(x;t)=h_{r}(\xi)$,
and therefore,
by Jacobi-Trudi identity \thetag{\ref{eq:JT}},
this implies $T_{\lambda}(x;t)=s_{\lambda}(\xi)$.
By \thetag{\ref{Schur-prod}} we obtain
\[
\sum_{\lambda}T_{\lambda}(x;q,t)s_{\lambda}(x)
=\prod_{i,j\geq1}
\frac{(tx_iy_j;q)_{\infty}}{(x_iy_j;q)_{\infty}}.
\]
The rest of the arguments is almost the same as in the proof of Corollary~\ref{conj:BigSchur}.
$\Box$
\end{demo}

\begin{remark}
Prof. A.~Lascoux said that Corollary~\ref{conj:BigSchur} and Corollary~\ref{conj:Macdonald}
are obtained as corollaries of Theorem~\ref{conj:general} by $\lambda$-ring arguments.
We cite his comment here.
What we need to do is to modify the argument of the symmetric functions.
For Corollary~\ref{conj:BigSchur},
we pass the argument $X$ to $X(1-t)= X-tX$.
For Corollary~\ref{conj:Macdonald},
we use $X(1-t)/(1-q)$.
This defines two transformations on the complete symmetric functions,
and therefore transformations on all the other functions.
In particular, for the power sums, it transforms
$p_k \to  (1-t^k)p_k $
or $p_k \to  (1-t^k)p_k/(1-q^k) $.
Any identity on symmetric functions which is valid for an infinite
alphabet $X$ remains valid for $X(1-t)$ and $X(1-t)/(1-q)$
and thus Theorem~\ref{conj:general} implies Corollary~\ref{conj:BigSchur} and Corollary~\ref{conj:Macdonald}.
About the $\lambda$-rings, the reader can consult \cite{La}.
\end{remark}


We also checked the Hall-Littlewood functions case,
and could not find a formula for the general case,
but found some nice formulas if we substitute $-1$ for $t$.

\begin{conjecture}
Let
\begin{equation*}
w(x;t)=\sum_{\lambda}\omega(\lambda)P_{\lambda}(x;t),
\end{equation*}
where $P_{\lambda}(x;t)$ denote the Hall-Littlewood function corresponding to the partition $\lambda$,
and the sum runs over all partitions $\lambda$.
Then
\begin{align*}
\log w(x;-1)+\sum_{n\geq1\text{ odd}}\frac1{2n}a^nc^np_{2n}
+\sum_{n\geq2\text{ even}}\frac1{2n}a^{\frac{n}2}c^{\frac{n}2}(a^{\frac{n}2}c^{\frac{n}2}-2b^{\frac{n}2}d^{\frac{n}2})p_{2n}\\
\in\Bbb{Q}[[p_1,p_3,p_5,\dots]]
\end{align*}
would hold.
$\Box$
\end{conjecture}
We might replace the Hall-Littlewood functions $P_{\lambda}(x;t)$
by the Macdonald polynomials $P_{\lambda}(x;q,t)$ in this conjecture.
Let $P_{\lambda}(x;q,t)$ denote the Macdonald polynomial corresponding to the partition $\lambda$
(See \cite{Ma}, IV, sec.4).
\begin{conjecture}
Let
\begin{equation*}
w(x;q,t)=\sum_{\lambda}\omega(\lambda)P_{\lambda}(x;q,t).
\end{equation*}
Here the sum runs over all partitions $\lambda$.
Then
\begin{align*}
\log w(x;q,-1)+\sum_{n\geq1\text{ odd}}\frac1{2n}a^nc^np_{2n}
+\sum_{n\geq2\text{ even}}\frac1{2n}a^{\frac{n}2}c^{\frac{n}2}(a^{\frac{n}2}c^{\frac{n}2}-2b^{\frac{n}2}d^{\frac{n}2})p_{2n}\\
\in\Bbb{Q}[[p_1,p_3,p_5,\dots]]
\end{align*}
would hold.
$\Box$
\end{conjecture}

\bigbreak
\noindent
{\bf Acknowledgment:}
The author would like to express his deep gratitude 
to Prof. R.~Stanley, Prof. S.~Okada and Prof. H.~Tagawa
for their helpful suggestions and discussions.

%
%

%
%

\medskip
\parindent=0mm

Masao ISHIKAWA, Department of Mathematics, Faculty of Education, \\
Tottori University, Tottori 680 8551, Japan

E-mail address: ishikawa@fed.tottori-u.ac.jp

\end{document}